\documentclass[11pt,a4paper]{article}
\usepackage[latin1]{inputenc}
\usepackage{amsmath,enumerate,fancyhdr,amssymb, amsthm, pstricks, epsf, lscape}
\usepackage{graphicx}

\newcommand{\tQ}{Q_R} 
\newcommand{\hQ}{Q_N} 

\newcommand{\lin}{\operatorname{lin}}
\newcommand{\gm}{G_{m, n}(M)}
\newcommand{\gmp}{G_{m, n}'(M)}

\newcommand{\width}{{\rm width}}

\newcommand{\eps}{\epsilon} 

\newcommand{\zad}[1]{\mathbb{Z}^{#1}}
\newcommand{\rad}[1]{\mathbb{R}^{#1}}

\def\eref#1{(\ref{#1})}

\newcommand{\commentout}[1]{}
\newcommand{\co}[1]{}
\co{
\setlength{\oddsidemargin}{1cm} 
\setlength{\evensidemargin}{1cm}
\setlength{\textwidth}{16.0cm} 
\setlength{\topmargin}{-1cm}
\setlength{\textheight}{21cm}
}
\co{
\setlength{\oddsidemargin}{0cm} 
\setlength{\evensidemargin}{0cm}
\setlength{\textwidth}{16.0cm} 
\setlength{\topmargin}{1cm}
\setlength{\textheight}{19cm}
}

\newtheorem{Definition}{Definition}

\newtheorem{Setup}{Setup}
\newtheorem{Example}{Example}
\newtheorem{Counterexample}{Counterexample}
\newtheorem{Proposition}{Proposition}
\newtheorem{Lemma}{Lemma}
\newtheorem{Theorem}{Theorem}
\newtheorem{Corollary}[Definition]{Corollary}
\newtheorem{Remark}[Definition]{Remark}
\newtheorem{Assumption}{Assumption}
\newtheorem{Recipe}{Recipe}

\newcommand{\beq}{\begin{equation}}
\newcommand{\eeq}{\end{equation}}
\newcommand{\beqa}{\begin{eqnarray}}
\newcommand{\eeqa}{\end{eqnarray}}
\newcommand{\ba}{\begin{array}}
\newcommand{\ea}{\end{array}}
\newcommand{\bac}{\begin{array}{ccccccccccc}}
\newcommand{\eac}{\end{array}}
\newcommand{\bprop}{\begin{Proposition}}
\newcommand{\eprop}{\end{Proposition}}

\newcommand{\bcex}{\begin{Counterexample}}
\newcommand{\ecex}{\end{Counterexample}}

\newcommand{\beqast}{\begin{eqnarray*}}
\newcommand{\eeqast}{\end{eqnarray*}}
\newcommand{\benum}{\begin{enumerate}}
\newcommand{\eenum}{\end{enumerate}}
\newcommand{\bit}{\begin{itemize}}
\newcommand{\eit}{\end{itemize}}
\newcommand{\bth}{\begin{Theorem}}
\newcommand{\enth}{\end{Theorem}}

\newcommand{\mydet}{\operatorname{det}}
\newcommand{\norm}[1]{\parallel \! #1 \! \parallel}

\newcommand{\latt}[1]{\mathbb{L}(#1)}

\newcommand{\bsetup}{\begin{Setup}}
\newcommand{\esetup}{\end{Setup}}
\newcommand{\ble}{\begin{Lemma}}
\newcommand{\ele}{\end{Lemma}}
\newcommand{\bex}{\begin{Example}}
\newcommand{\eex}{\end{Example}}
\newcommand{\bcor}{\begin{Corollary}}
\newcommand{\ecor}{\end{Corollary}}
\newcommand{\brem}{\begin{Remark}}
\newcommand{\erem}{\end{Remark}}
\newcommand{\bass}{\begin{Assumption}}
\newcommand{\eass}{\end{Assumption}}

\newcommand{\brep}{\begin{Recipe}}
\newcommand{\erep}{\end{Recipe}}

\newcommand{\la}{\langle}
\newcommand{\ra}{\rangle}

\newcommand{\nin}{\noindent}

\newcommand{\pf}[1]{\vspace{.35cm} \nin {\bf Proof {#1} }}

\newcommand{\bpx}{\begin{pmatrix}}
\newcommand{\epx}{\end{pmatrix}}

\newcommand{\bbx}{\begin{bmatrix}}
\newcommand{\ebx}{\end{bmatrix}}

\newcommand{\bb }{B\&B \ } 

\co{

math: CO (combinatorics); OC (optimization and control); 


ACM   F.2.1 Numerical Algorithms and Problems   (includes Number-theoretic computations (e.g., factoring, primality testing) 
  F.2.2 Nonnumerical Algorithms and Problems  (includes   Computations on discrete structures and Geometrical problems and computations) 
}

\begin{document}
\date{}
\title{Basis Reduction, and the Complexity of Branch-and-Bound}
\author{G\'{a}bor Pataki and Mustafa Tural \thanks{Department of Statistics and Operations Research, UNC Chapel Hill, {\bf gabor@unc.edu, tural@email.unc.edu}} \\
Department of Statistics and Operations Research, UNC Chapel Hill}

\maketitle

\thispagestyle{empty}
\co{
The classical branch-and-bound algorithm for the integer feasibility problem
\begin{equation} \label{ip} 
\text{Find} \, x \in Q \cap \zad{n}, \,\, \text{with} \,\, Q = \left\{ \, x  \, | \, \begin{pmatrix} \ell_1 \\ \ell_2 \end{pmatrix} \, \leq \, \begin{pmatrix} A \\ I \end{pmatrix}  \leq \begin{pmatrix} w_1 \\ w_2 \end{pmatrix} \right\}
\eeq
has exponential worst case complexity. We prove that it is  surprisingly efficient on reformulations of (\ref{ip}), in which the columns of the constraint matrix are short, and near orthogonal, i.e. a reduced basis of the generated lattice.

The analysis is in the spirit of Furst and Kannan's work on the subset sum problem, and it uses their ideas to bound the number of integral matrices 
for which the shortest nonzero vectors of certain lattices are long. We also 
use an upper bound on the size of the branch-and-bound tree. 

We also explore practical aspects of this result.
We compute numerical values of $M$ which guarantee that $90$, and $99$ percent of the reformulated problems solve at the root: these turn out to be 
surprisingly small when the problem size is moderate. A computational study also confirms that random integer programs become  easier, as the coefficients grow. 
}

\vspace{-0.6cm}

\begin{abstract}

The classical branch-and-bound algorithm for the integer feasibility problem
\begin{equation} \label{ip0} 
\text{Find} \, x \in Q \cap \zad{n}, \,\, \text{with} \,\, Q = \left\{ \, x  \, | \, \begin{pmatrix} \ell_1 \\ \ell_2 \end{pmatrix} \, \leq \, \begin{pmatrix} A \\ I \end{pmatrix} x \leq \begin{pmatrix} w_1 \\ w_2 \end{pmatrix} \right\}
\eeq
has exponential worst case complexity. 

We prove that it is  surprisingly efficient on reformulations of (\ref{ip0}), 
in which the columns of the constraint matrix are short, and near orthogonal, i.e. a reduced basis of the generated lattice; 
when the entries of $A$ are from $\{1, \dots, M\}$  for a large enough $M$, branch-and-bound solves almost all reformulated instances at the rootnode.
For all $A$ matrices we prove an upper bound on the width of the reformulations along the last unit vector.

The analysis builds on the ideas of Furst and Kannan to bound the number of integral matrices 
for which the shortest nonzero vectors of certain lattices are long, and also uses
a bound on the size of the branch-and-bound tree based on the norms of the Gram-Schmidt vectors 
of the constraint matrix. 

We explore practical aspects of these results.
First, we compute numerical values of $M$ which guarantee that $90$, and $99$ percent of the reformulated problems solve at the root: these turn out to be 
surprisingly small when the problem size is moderate. 
Second, we confirm with a computational study that random integer programs become easier, as the coefficients grow.

\end{abstract}

\vspace{-0.3cm}

\section{Introduction and main results} 

\co{Given a polyhedron $Q$ described by inequalities, the Integer Programming (IP) feasibility problem asks 
whether $Q$ contains an integral point. 
Branch-and-bound, which we abbreviate as  B\&B is a classical solution method. 
It starts with $Q$ as the sole subproblem (node). In a general step, one chooses a subproblem $Q'$, a variable $x_i$, 
and creates nodes $Q'\cap\left\lbrace x | x_i=\gamma\right\rbrace $, where $\gamma$ ranges over all possible integer values of  $x_i$.
We  repeat this until all subproblems are found to be empty, or an integral point is found in one of them.
}

The Integer Programming (IP) feasibility problem asks whether a polyhedron $Q$ contains an integral point. 
Branch-and-bound, which we abbreviate as  B\&B is a classical solution method, first proposed by Land and Doig in \cite{LD60}.
It starts with $Q$ as the sole subproblem (node). In a general step, one chooses a subproblem $Q'$, a variable $x_i$, 
and creates nodes $Q'\cap\left\lbrace x | x_i=\gamma\right\rbrace $, where $\gamma$ ranges over all possible integer values of  $x_i$.
We  repeat this until all subproblems are shown to be empty, or we find an integral point in one of them.

B\&B (and its version used to solve optimization problems) enhanced by cutting planes is a dependable algorithm implemented in most commercial software packages.
However, instances in \cite{J74, C80, GNS98Complex, KP09, AL04, AL06} 
show that it is theoretically inefficient: it can take an exponential number of 
subproblems to prove the infeasibility of simple knapsack problems. While B\&B is inefficient in the worst case, Cornu{\'{e}}jols et al. in \cite{CKL06} developed
useful computational tools to give an early estimate on the size of the B\&B tree in practice.

Since IP feasibility is NP-complete, one can ask for polynomiality of a solution method only in fixed dimension. 
All algorithms that achieve such complexity rely on advanced techniques. 
The algorithms of Lenstra \cite{L83} and Kannan \cite{K87} first round the polyhedron (i.e. apply a transformation to 
make it have a spherical appearance), then use basis reduction to reduce the problem 
to a provably small number of smaller dimensional subproblems. On the subproblems the algorithms are applied recursively, e.g. rounding is done again. 
Generalized basis reduction, proposed by Lov\'asz and Scarf in \cite{LS92}
avoids rounding, 
but needs to solve a sequence of linear programs to create the subproblems. 
\co{In fixed dimension 
one can even {\em count} the 
number of feasible solutions in polynomial time: see the papers of Barvinok \cite{Bar94}, Dyer and Kannan \cite{DyKa97}, 
and Koeppe \cite{Koep07}. 
We refer to \cite{CRSS93, DHTY04} for successful implementations of these theoretically efficient methods. }
\co{, and to 
Haus et al. \cite{HKW03} for a finite augmentation type algorithm using basis reduction. }

\co{We refer to \cite{CRSS93} and De Loera et al. 
\cite{DHTY04} for successful implementations of these theoretically efficient methods, and to 
Haus et al \cite{HKW03} to a finite augmentation type algorithm using basis reduction. }

\co{
A simpler use of basis reduction in integer programming is preprocessing \eref{ip0} to create an instance with 
short and near orthogonal columns in the constraint matrix, then simply feeding it to an IP solver.
We describe two such methods that were  proposed recently. }

There is a simpler way to use basis reduction in integer programming: preprocessing \eref{ip0} to create an instance with 
short and near orthogonal columns in the constraint matrix, then simply feeding it to an IP solver.
We describe two such methods that were  proposed recently. 
\co{
\beq \label{ip} 
\text{Find} \, x \in Q \cap \zad{n}, \,\, \text{with} \,\, Q = \left\{ \, x  \, | \, \bpx \ell_1 \\ \ell_2 \epx \, \leq \, \bpx A \\ I \epx x \leq \bpx w_1 \\ w_2 \epx \right\},
\eeq}
We assume that $A$ is an integral matrix with $m$ rows, and $n$ columns, and the $w_i$ and $\ell_i$ are integral vectors.

The rangespace reformulation of \eref{ip0} proposed by Krishnamoorthy and Pataki in \cite{KP09} is
\beq \label{ipt} 
\text{Find} \, y \in Q_R \cap \zad{n}, \,\, \text{with} \,\, \tQ = \left\{ \, y  \, | \, \bpx \ell_1 \\ \ell_2 \epx \, \leq \, \bpx A \\ I \epx Uy \leq \bpx w_1 \\ w_2 \epx \right\},
\eeq
where $U$ is a unimodular matrix computed to make the columns of the constraint matrix a reduced basis of the generated lattice. 

The nullspace reformulation of Aardal, Hurkens, and Lenstra proposed in \cite{AHL00}, and further studied in \cite{ABHLS00}
is applicable, when the rows of $A$  are linearly independent, and $w_1 = \ell_1$. 
It is 
\beq \label{iptt} 
\text{Find} \, y \in Q_N \cap \zad{n-m}, \,\, \text{with} \,\, \hQ = \left\{ \, y  \, | \, \ell_2 - x_0 \, \leq \, By \, \leq \, w_2 - x_0 \right\},
\eeq
where $x_0 \in \zad{n}$ satisfies $Ax_0 = \ell_1, \,$ and the columns of $B$ are a reduced basis of the lattice
\mbox{$\{ \, x \in \zad{n} \, | \, Ax = 0 \, \}$}. 

We analyze the use of
Lenstra-Lenstra-Lov\'asz (LLL) \cite{LLL82}, 
and reciprocal Korkhine-Zolotarev (RKZ) reduced bases  \cite{LLS90} in the reformulations, and use 
Korkhine-Zolotarev (KZ) reduced bases \cite{K87}, \cite{KZ1873} in our computational study. 
The definitions of these reducedness concepts are given in Section \ref{proofs}. 
When $Q_R$ is computed using LLL reduction, we call it the LLL-rangespace reformulation of $Q$, and abusing notation we also call 
\eref{ipt} the LLL-rangespace reformulation of \eref{ip0}. 
Similarly we talk about 
LLL-nullspace, RKZ-rangespace, and RKZ-nullspace reformulations. 

\bex \label{ex1}
The polyhedron 
\beq \label{ex1-prob}
\ba{rcl}
207 \, \leq & 41 x_1 + 38 x_2 & \leq \, 217 \\
0   \, \leq & x_1, x_2 & \leq \, 10 \\
\ea
\eeq
is shown on the first picture of Figure \ref{fig1}. It is long and thin, and defines an infeasible, and 
relatively difficult integer feasibility problem for B\&B, 
as branching on either $x_1$ or $x_2$ yields $6$ subproblems. 
Lenstra's and Kannan's algorithms would first transform this polyhedron to make it more spherical;
generalized basis reduction  would solve a sequence of linear programs to find the direction $x_1 + x_2$ along which the polyhedron is thin.

The LLL-rangespace reformulation is 
\beq
\ba{rcl}
207 \, \leq & -3 x_1 + 8 x_2 & \leq \, 217 \\
0 \, \leq & - x_1 - 10 x_2 & \leq \, 10 \\
0 \, \leq & x_1 + 11 x_2 & \leq \, 10
\ea
\eeq
shown  on the second picture of Figure \ref{fig1}: now branching on $y_2$ proves integer infeasibility. (A similar example was given in 
\cite{KP09}). 
\begin{center}
\begin{figure}[ht]
\includegraphics[scale=0.55, angle=-90, clip=true, trim=2.4cm 2.4cm 8.5cm 0]{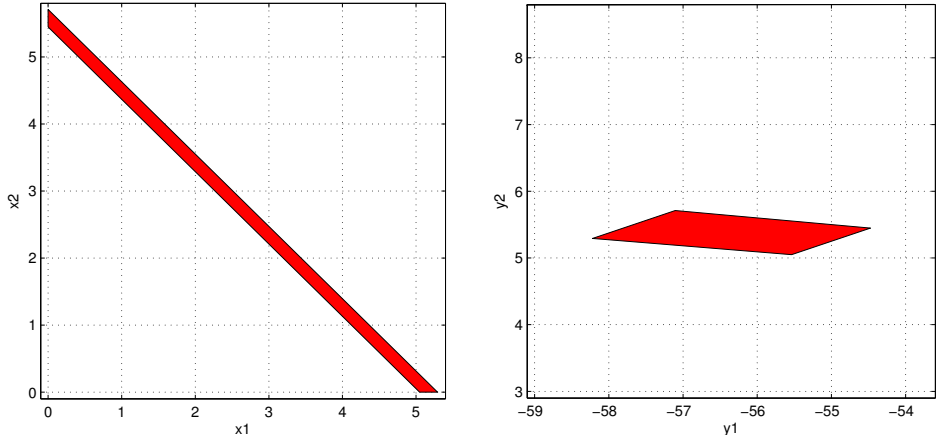}
\caption{The polyhedron of Example \ref{ex1} before and after the reformulation} 
\label{fig1}
\end{figure}
\end{center}
\vspace{-1.2cm}
\co{
\begin{center}
\begin{figure}[ht]
\label{bla} 
\includegraphics[scale=0.6, angle=-90, clip=true, trim=1.5cm 0 7.5cm 0]{comb2.eps}
\end{figure}
\end{center}
}
\eex
\vspace{-0.4cm}
The reformulation methods are easier to describe, than, say Lenstra's algorithm, and 
are also successful in practice in solving several classes of hard integer programs: see \cite{AHL00, ABHLS00, KP09}.
For instance, the original formulations of the marketshare problems of Cornu{\'{e}}jols and Dawande in \cite{CD99}
are notoriously difficult for commercial solvers, while 
the nullspace reformulations are much easier to solve as shown by Aardal et al in \cite{ABHLS00}.

However, they seem difficult to analyze in general. 
Aardal and Lenstra in \cite{AL04, AL06} studied knapsack problems with a nonnegativity constraint, and 
the constraint vector $a$ having a given decomposition 
$a=\lambda p + r, \,$ with $p$ and $r$ integral vectors, and $\lambda$ an integer, large compared to $\norm{p}$ and $\norm{r}$. 
They proved a lower bound on the norm of the last vector in the nullspace reformulation, and argued that
branching on the corresponding variable will create a small number of B\&B nodes. 
Krishnamoorthy and Pataki in \cite{KP09} pointed out a gap in this argument, 
and showed that branching on the constraint $px$ in $Q$ (which creates a small number of subproblems, as $\lambda$ is large), 
is equivalent to branching on the last variable in $\tQ$ and $\hQ$. 
\co{In \cite{PT07} the authors of this paper analyzed knapsack problems with the constraint vector satisfying 
$\norm{a} \geq 2^{n^2/2}$, but assuming no a priori structure. They proved a similar result, as 
in \cite{KP09}, i.e.  that a small number of subproblems are created by branching on the last variable in the reformulations. }
\co{
They showed that from the inverse of $U$ in the rangespace reformulation one can extract 
a vector $p$ which is near parallel to $a, $ leading again to a small number of subproblems
when branching on the last variable. }

\co{
The only analysis that exists so far is for knapsack problems with a constraint vector of the form 
$a=\lambda p + r, \,$ with $p$ and $r$ integral vectors, and $\lambda$ an integer, large compared to $\norm{p}$ and $\norm{r}$. 
Aardal and Lenstra in \cite{AL04, AL06} proved a lower bound on the norm of the last vector in the nullspace reformulation, and argued that 
branching on such a long vector creates a small number of B\&B nodes. 
}

\co{

Aa
analyzed the nullspace reformulation method on an equality constrained knapsack problem 
\begin{align}
(\lambda p+r) x & = \beta \\
x & \geq 0,
\end{align}
where 

They proved a lower bound on the 
}
\co{
Following Aardal and Lenstra in \cite{AL04, AL06},  Krishnamoorthy and Pataki in \cite{KP09} 
analysed knapsack problems with a constraint vector $a$ having a given decomposition 
$a=\lambda p + r, \,$ with $p$ and $r$ integral vectors, and $\lambda$ an integer, large compared to $\norm{p}$ and $\norm{r}$. 
They show that branching on the constraint $px$ in $Q$ (which creates a small number of subproblems, by virtue of $\lambda$ being large) 
is equivalent to branching on the last variable in $\tQ$ and $\hQ$. 
}
A result one may hope for is 
proving polynomiality of B\&B  on the reformulations of \eref{ip0} when the dimension is fixed. 
While this seems difficult, 
we give a different, and perhaps even more surprising  complexity analysis. It is 
in the spirit of Furst and Kannan's work in \cite{FK89} on subset sum problems and builds on a generalization of their Lemma 1 
to bound the fraction of integral matrices for which the shortest nonzero vectors of certain corresponding lattices are short.
We also use 
an upper bound on the size of the B\&B tree, 
which depends on the norms of the Gram-Schmidt vectors 
of the constraint matrix. 
We introduce necessary notation, and state our results, then give a comparison with \cite{FK89}.

When a statement is true for all, but at most a fraction of $1/2^n$ of the elements of a set $S$, we say that it is true for {\em almost all} elements.
The value of $n$ will be clear from the  context. 
{\em Reverse  B\&B } is B\&B branching on the variables in reverse order starting with the one of highest index. 
We assume $w_2 > \ell_2$, and for simplicity of stating the results we also assume $n \geq 5.$ For positive integers $m, \, n$ and $M$ 
we denote by $\gm$ the set of matrices with $m$ rows, and $n$ columns, and the entries drawn from 
$\{1, \dots, M \}$. We denote by $\gmp$ the subset of $\gm$ consisting of matrices with linearly independent rows, and let 
\begin{align}
\chi_{m,n}(M) & = \dfrac{|\gmp|}{|\gm|}.
\end{align}
It is shown by Bourgain et. al in \cite{BVW09} that $\chi_{m,m}(M)$ (and therefore also 
$\chi_{m,n}(M)$ for $m \leq n$) are of the order $1 - o(1)$. In this paper we will use $\chi(m,n,M) \geq 1/2$ for simplicity. 

For matrices (and vectors) $A$ and $B$, we write $(A; B)$ for $\bpx A \\ B \epx$.
For an $m$ by $n$ integral matrix $A$ with independent rows we write $\gcd(A)$ for the greatest common divisor of the $m$ by $m$ subdeterminants of $A$. 
If B\&B generates at most one node at each level of the tree, we say that it solves an integer feasibility problem at the rootnode.

If $Q$ is a polyhedron,  and $z$ is an integral vector, then the width of $Q$ along $z$ is
\begin{align}
\width(z,Q) & = \max \left\lbrace \, \la z, x \ra \, | \, x \in Q \right\rbrace - \min \left\lbrace \, \la z, x \ra \, | \, x \in Q \right\rbrace. 
\end{align}
The main results of the paper follow.
\bth \label{thm-rkz} There are positive constants $d_1 \leq 2, \,$ and $d_2 \leq 12 \,$ such that the following hold.
\benum 
\item \label{thm-rkz-r} If 
\beq
M > (d_1 n \norm{(w_1; w_2) - (\ell_1; \ell_2)})^{n/m+1},
\eeq
then for almost all $A \in \gm$ reverse B\&B solves the RKZ-rangespace reformulation of \eref{ip0} at the rootnode.
\item \label{thm-rkz-n} If 
\beq
M > (d_2 (n-m) \norm{ w_2 - \ell_2})^{n/m},
\eeq
then for almost all $A \in \gmp$ reverse B\&B solves the RKZ-nullspace reformulation of \eref{ip0} at the rootnode. \qed
\eenum
\enth

The proofs also show that when $M$ obeys the above bounds, then $Q$ has at most one element for almost all $A \in \gm$. 
When $n/m$ is fixed, and the problems are binary, the magnitude of $M$ required is a polynomial in $n$.

\bth \label{thm-lll} The conclusions of Theorem \ref{thm-rkz} hold for the LLL-reformulations, if the bounds on $M$ are 
\beq \nonumber
(2^{(n+4)/2} \norm{(w_1; w_2) - (\ell_1; \ell_2)})^{n/m+1}, \; \text{and} \; (2^{(n-m+4)/2} \norm{w_2 - \ell_2})^{n/m},
\eeq
respectively. 
\qed
\enth

Furst and Kannan, based on Lagarias' and Odlyzko's \cite{LO85} and Frieze's  \cite{F86} work
show that the subset sum problem is solvable in polynomial time using a simple iterative method 
for almost all weight vectors in $\{ \, 1,\dots,M \}^n,$  and all right hand sides, 
when $M$ is sufficiently large, and a reduced basis of the orthogonal lattice of the weight vector is  available. 
The lower bound on $M$ is $2^{c n \log n}, \,$ when the basis is RKZ reduced, and 
$2^{d n^2}, \,$ when it is LLL reduced. Here $c$ and $d$ are positive constants.

Theorems \ref{thm-rkz} and \ref{thm-lll}   generalize the solvability results from subset sum problems to bounded integer programs; also, 
we prove them via branch-and-bound, an algorithm  considered inefficient from the theoretical point of view. 

Proposition \ref{widthcor} gives another indication why the reformulations are relatively easy. One can observe that 
$\mydet (A A^T)$ can be quite large even for moderate values of $M$, 
if $A \in \gm$ is a random matrix with $m \leq n$, although we could not find any theoretical studies on 
the subject.  For instance, for a random $A \in G_{4,30}(100)$ we  found $\mydet (A A^T)$ to be of the order $10^{18}. \,$ 

While we cannot give a tight upper bound 
on the size of the B\&B tree in terms of this determinant, we are able to bound the width of the reformulations along the last unit vector
for any $A$ (i.e. not just almost all). 

\bprop \label{widthcor} 
If $Q_R$ and $Q_N$ are computed using RKZ reduction, then
\begin{align} \label{widthrangerkz} 
\width(e_n, Q_R) & \leq  \dfrac{\sqrt{n} \norm{(w_1; w_2) - (\ell_1; \ell_2)}}{\mydet (A A^T+I)^{1/(2n)}}.
\end{align}
Also, if $A$ has independent rows, then 
\begin{align}
\label{widthnullrkz} 
\width(e_{n-m}, Q_N) & \leq \dfrac{\gcd(A) \sqrt{n-m} \norm{ w_2 - \ell_2}}{\mydet (AA^T)^{1/(2(n-m))}}.
\end{align}
\eprop

The same results hold for the LLL-reformulations, if 
$\sqrt{n}$ and $\sqrt{n-m}$ are replaced by $2^{(n-1)/4}$ and $2^{(n-m-1)/4}$, respectively. 
\qed

\brem
\co{By the correspondence of the widths between branching directions, } As described in Section 5 of \cite{KP09}, and in \cite{ML04} for the nullspace reformulation,
directions achieving the same widths exist in $Q, \,$ and they can be quickly computed. 
For instance, if $p$ is the last row of $U^{-1}, \,$ then 
$\width(e_n, Q_R) = \width(p, Q)$. 
\erem

\co{
\bcor \label{widthcor2} 
Suppose $n$ is fixed. Then we 
can compute in polynomial time $f_1 \in \zad{n}$ such that 
\begin{align} \label{widthrangerkz2} 
\width(f_1, Q) & \leq  \dfrac{\sqrt{n} \norm{(w_1; w_2) - (\ell_1; \ell_2)}}{\mydet (A A^T+I)^{1/(2n)}}.
\end{align}
If $w_1 = \ell_1, \,$ then we can compute in polynomial time $f_2 \in \zad{n}$

RKZ reduced basis of the lattice $\{ \, x \, | \, Ax = 0 \, \}, \,$ 
generated by the constraint matrix of $Q$ we can find $f_1 \in \zad{n}$ with

When $n$ is fixed, we can compute in poly
We can compute There exists $d \in \zad{n}$ such that 
If $Q_R$ and $Q_N$ are computed using RKZ reduction, then
\begin{align} \label{widthrangerkz} 
\width(e_n, Q_R) & \leq  \dfrac{\sqrt{n} \norm{(w_1; w_2) - (\ell_1; \ell_2)}}{\mydet (A A^T+I)^{1/(2n)}}.
\end{align}
Also, if $A$ has independent rows, then 
\begin{align}
\label{widthnullrkz} 
\width(e_{n-m}, Q_N) & \leq \dfrac{\gcd(A) \sqrt{n-m} \norm{ w_2 - \ell_2}}{\mydet (AA^T)^{1/(2n)}}.
\end{align}
\eprop

The same results hold for the LLL-reformulations, if 
$\sqrt{n}$ and $\sqrt{n-m}$ are replaced by $2^{(n-1)/4}$ and $2^{(n-m-1)/4}$, respectively. 
\qed
}

A practitioner of integer programming may ask for the value of Theorems \ref{thm-rkz} and \ref{thm-lll}.
Proposition \ref{actual} 
and a computational study put these results into a more practical perspective.
Proposition \ref{actual} shows that when $m$ and $n$ are not too large,
already fairly small values of $M$ guarantee that 
the RKZ-nullspace reformulation (which has the smallest bound on $M$) 
of the majority of binary integer programs get solved at the rootnode. 

\bprop \label{actual}
Suppose that $m$ and $n$ are chosen according to 
Table \ref{mnMtable90}, and $M$ is as shown in the third column. 
\begin{table}[h]
\begin{center}
\begin{tabular}{ | c | c | c | c | p{5cm} |}
\hline
$n$ & $m$ &  $M$ for $90$ \% & $M$ for $99$ \% \\ \hline
$30$ & $20$ &  $33$          &    $37$      \\ \hline
$50$ & $20$ & $1912$         &   $2145$    \\ \hline
$50$ & $30$ & $96$           & $103$       \\ \hline
$60$ & $30$ & $420$          &  $454$      \\ \hline
$70$ & $40$ & $197$          &  $209$      \\ \hline
\end{tabular}
\end{center}
\caption{Values of $M$ to make sure that the RKZ-nullspace reformulation of $90$ or $99$ $\%$ of the instances of type \eref{Ab} solve at the rootnode} 
\label{mnMtable90} 
\end{table}
Then for at least $90 \%$ of \mbox{$A \in \gmp, \,$} and all $b$ right hand sides, 
reverse B\&B solves the RKZ-nullspace reformulation of 
\beq \label{Ab} 
\ba{rcl}
Ax & = & b \\
x & \in & \{ 0, 1 \}^n 
\ea
\eeq
at the rootnode. The same is true for $99 \%$ of $A \in \gmp, \,$ if 
$M$ is as shown in the fourth column.  \qed
\eprop
\vspace{-0.2cm}
Note that $2^{n-m}$ is the best upper bound one can give on the 
number of nodes when B\&B is run on the original formulation \eref{Ab}; also, randomly generated IPs with $n-m=30$ are nontrivial 
 even for commercial solvers. 

\co{

in \eref{widthrangerkz} is replaced by ; and 
\begin{align}
\width(e_n, Q_R) & \leq \dfrac{2^{(n-1)/4} \norm{(w_1, w_2) - (\ell_1, \ell_2)}}{\mydet (A A^T+I)} \\
\width(e_n, Q_R) & \leq  \dfrac{2^{(n-m-1)/4} \norm{ w_2 - \ell_2}}{\mydet A A^T}. 
\end{align}
\eprop
}

According to Theorems \ref{thm-rkz} and \ref{thm-lll}, random integer programs with coefficients drawn from $\{\, 1, \dots, M \, \}$
should get easier, as $M$ grows. Our computational study confirms this somewhat counterintuitive 
hypothesis on the family of marketshare problems of Cornu{\'{e}}jols and Dawande  in \cite{CD99}.

We generated twelve $5$ by $40$ matrices with entries drawn from $\{1, \dots, M \}$ with 
$M = 100, 1000, \,$ and $10000$ (this is $36$ matrices overall), set $b = \lfloor Ae/2 \rfloor, \,$ where $e$ is the vector of all ones, and constructed the instances 
of type \eref{Ab}, and 
\beq \label{Abineq} 
\ba{rcl}
b - e  \leq  & Ax  & \leq  b \\
      & x  \in  \{ 0, 1 \}^n. & 
\ea
\eeq
The latter of these are a relaxed version, which correspond to trying to find an almost-equal market split. 

Table \ref{bbtable} shows the average number of nodes that the commercial IP solver CPLEX 9.0 took
to solve the rangespace reformulation 
of the inequality- and the nullspace reformulation of the equality constrained problems. 

\begin{table}[h]
\begin{center}
\begin{tabular}{| l | r | r | }
\hline
$M$   & EQUALITY &  INEQUALITY \\ \hline
$100$ & $17531.92$ &  $38884.92$ \\ \hline
$1000$ & $1254.42$ & $22899.67$ \\ \hline
$10000$ & $200.83$ & $1975.67$ \\ \hline
\end{tabular}
\end{center}
\caption{Average number of \bb nodes to solve the inequality- and equality-constrained marketshare problems} 
\label{bbtable} 
\end{table}

Since RKZ reformulation is not implemented in any software that we know of, 
we used the Korkhine-Zolotarev (KZ) reduction routine from the NTL library \cite{NTL}. 
For brevity we only report the number of B\&B nodes, and not the actual computing times. 

All equality constrained instances turned out to be infeasible, except two, corresponding to $M=100$. 
Among the inequality constrained problems there were 
fifteen  feasible ones: all twelve  with $M=100, \,$ and three with $M=1000$. Since infeasible problems tend to be harder, this explains 
the more moderate decrease in difficulty as we go from $M=100$ to $M=1000$. 

Table \ref{bbtable} confirms the theoretical findings of the paper: 
the reformulations of random integer programs become easier as the size of the coefficients grow.

In Section \ref{proofs} we introduce further necessary notation, and 
give the proof of Theorems \ref{thm-rkz} and \ref{thm-lll}.

\vspace{-0.3cm}

\section{Further notation, and proofs} 
\label{proofs} 

A lattice is a set of the form 
\beq \label{def-latt-B}
L \, = \, \latt{B} \, = \, \{ \, Bx \, | \, x \in \mathbb{Z}^{r} \, \},
\eeq
where $B$ is a real matrix with $r$ independent columns, called a {\em basis} of $L$, 
and $r$ is called the {\em rank} of $L$. 

The euclidean norm of a shortest nonzero vector in $L$ is denoted by $\lambda_1(L)$, and Hermite's constant is 
\beq 
C_j=\sup\left\lbrace \lambda_1(L)^2/ (\mydet L)^{2/j} \, | \, L \text{ is a lattice of rank $j$ }\right\rbrace.
\eeq
We define 
\beq
\gamma_i = \max\left\lbrace C_1, \dots, C_i \right\rbrace.
\eeq
\co{
\beq 
C_j=\sup\left\lbrace \lambda(L)^2/ (\mydet L)^{2/j} \, | \, L \text{ is a lattice of rank $j$ }\right\rbrace.
\eeq

It is straightforward to prove that $C_j \leq 2j/3$, see for example \cite{LLS90}. A better upper bound on $C_j$ was given by Blichfeldt in [?]: ****

\begin{equation} \label{Bl}
C_j \leq \dfrac{2}{\pi} {\Gamma \left( {\dfrac{n+4}{2}} \right) }^{2/n}.
\end{equation}
}
A matrix $A \,$ defines two lattices that we are interested in: 
\beq
L_R(A) = \mathbb{L}(A;I),  \co{ {\bpx A \\ I \epx},} \, L_N(A) =\left\lbrace x \in \mathbb{Z}^{n} | Ax = 0 \right\rbrace,
\eeq
where we recall that $(A; I)$ is the matrix obtained by stacking $A$ on top of $I$. 

Given independent vectors $b_1, \dots, b_r $, the vectors $b_1^*, \dots, b_r^*$ form the Gram-Schmidt orthogonalization of 
$b_1, \dots, b_r, \,$ if $b_1^* = b_1, \,$ and $b_i^*$ is the projection of $b_i$ onto the orthogonal complement of the subspace spanned 
by $b_1, \dots, b_{i-1}$ for $i \geq 2$. We have 
\beq \label{bibist}
b_i =  b_i^* + \sum_{j=1}^{i-1} \mu_{ij} b_j^*,
\eeq
with 
\beq \label{muij} 
\ba{rcll} 
\mu_{ij} & = & \la b_i, b_j^* \ra /\norm{b_j^*}^2 & \,\,\, (1 \leq j < i \leq r).
\ea
\eeq
We call $b_1, \dots, b_r$ {\em LLL-reduced} if 
\beqa \label{mucond}
| \mu_{ij} | & \leq & \frac{1}{2} \; (1 \leq j < i \leq r), \\ \label{exch-cond}
\norm{\mu_{i, i-1} b_{i-1}^* + b_i^*}^2 & \geq & \frac{3}{4} \norm{b_{i-1}^*}^2  (1 < i \leq r).
\eeqa
An LLL-reduced basis can be computed in polynomial time for varying $n$.

Let 
\beq
b_i(k) =  b_i^* + \sum_{j=k}^{i-1} \mu_{ij} b_j^* \,\, (1 \leq k \leq i \leq r),
\eeq
and for $i=1, \dots, r \,$ let $L_i$ be the lattice generated by 
$$
b_i(i), b_{i+1}(i), \dots, b_r(i).
$$
We call $b_1, \dots, b_r$ {\em Korkhine-Zolotarev reduced} (KZ-reduced for short)
if $b_i(i)$ is the shortest nonzero vector in $L_i \,$ for all $i$. Since $L_1 = L \,$ and $\, b_1(1)=b_1, \,$
in a KZ-reduced basis the first vector is the shortest nonzero vector of $L$. Computing the shortest nonzero vector in a 
lattice is expected to be hard, though it is not known to be NP-hard. It can be done 
in polynomial time when the dimension is fixed, and so can be computing a KZ reduced basis.

Given a lattice $L$ its reciprocal lattice $L'$ is defined as 
$$
L' = \{ \, z \in \lin L \, | \, \la z, x \ra \in \zad{} \, \forall x \in L \, \}.
$$
For a basis $b_1, \dots, b_r$ of the lattice $L, \,$ there is a unique basis $b_1', \dots, b_r'$ of $L'$ called the {\em reciprocal basis} of 
$b_1, \dots, b_r$, with 
\begin{align} 
\la b_i, b_j' \ra & =  \begin{array}{rl} 1 & \text{if} \,\, i + j = r+1 \\
                               0 & \, \text{otherwise}.
                           \end{array} 
\end{align} 
We call a basis $b_1, \dots, b_r$ a reciprocal Korkhine Zolotarev (RKZ) basis of $L, \,$ if 
its reciprocal basis is a KZ reduced basis of $L'$. 
Below we collect the important properties of RKZ and LLL reduced bases.
\ble
Suppose that $b_1, \dots, b_r$ is a basis of the lattice $L \,$ with Gram-Schmidt orthogonalization 
$b_1^*, \dots, b_r^*. \,$ Then 
\benum
\item if $b_1, \dots, b_r$ is RKZ reduced, then 
\beq \label{bistci}
\norm{b_i^*} \geq \lambda_1(L)/C_i,
\eeq
and 
\beq \label{brrkz} 
\norm{b_r^*} \geq (\mydet L)^{1/r}/\sqrt{r}.
\eeq
\item if $b_1, \dots, b_r$ is LLL reduced, then 
\beq \label{bist2i}
\norm{b_i^*} \geq \lambda_1(L)/2^{(i-1)/2},
\eeq
and 
\beq \label{brlll} 
\norm{b_r^*} \geq (\mydet L)^{1/r}/2^{(r-1)/4}.
\eeq
\eenum
\ele
\pf{} Statement \eref{bistci} is proven in \cite{LLS90}. Let 
$b_1', \dots, b_r'$ be the reciprocal basis. Since $b_1'$ is the shortest nonzero vector of $L', \,$ Minkowski's theorem implies
\beq
\norm{b_1'} \leq \sqrt{r} (\mydet L')^{1/r}.
\eeq
Combining this with  $\norm{b_1'} = 1/\norm{b_r^*}, $ and $\mydet L' = 1/\mydet L$ prove \eref{brrkz}. 
Statement \eref{bist2i} was proven in \cite{LLL82}. Multiplying the inequalities 
\beq
\norm{b_i^*} \leq 2^{(r-i)/2} \norm{b_r^*} \; (i=1, \dots, r),
\eeq
and using $\norm{b_1^*} \dots \norm{b_r^*} = \mydet L$ gives \eref{brlll}. 
\qed

\ble \label{widthlemma} 
Let $P$ be a polyhedron
\beq \label{pol}
P = \left\lbrace y \in \rad{r} \, | \, \ell \leq By \leq w \right\rbrace,
\eeq
and $b_1^*, \dots, b_r^*$ the Gram-Schmidt orthogonalization of the columns of $B$.
When reverse B\&B is applied to $P$, the number of nodes on the level of $y_i$ is at most 
\beq
\left( \left\lfloor {\dfrac{\norm{w-\ell} }{\norm{b_i^*}}  }\right\rfloor +1\right) \dots \left( \left\lfloor \dfrac{\norm{w-\ell}}{\norm{b_r^*}}  \right\rfloor +1\right).
\eeq
\ele
\vspace{-0.4cm}
\pf{} 
First we show 
\beq \label{erwidth} 
\width(e_r, P) \leq \norm{w-\ell}/\norm{b_r^*}.
\eeq
\co{
\beq \label{erwidth} 
\norm{w-\ell}/\norm{b_r^*}.
\eeq}
Let $x_{r,1}$ and $x_{r,2}$ denote the maximum and the minimum of $x_r$ over $P$.
Writing 
$\bar{B}$ for the matrix composed of the first $r-1$ columns of $B$, and $b_r$ for the last column, it holds that 
there is $x_1, x_2 \in \mathbb{R}^{r-1}$ such that $\bar{B} x_1 + b_r x_{r,1} \,$ and $\bar{B} x_2 + b_r x_{r,2} \,$ are in $P$. So 
\co{\begin{align} 
\norm{w-\ell} & \geq \, \norm{B \left[ (x_1;x_{r,1})-(x_2;x_{r,2}) \right]} \, \geq \, \norm{b_r^*} \vert x_{r,1}-x_{r,2} \vert = \, \norm{b_r^*} \width(e_r,P)
}
\begin{align} \nonumber
\hspace{-0.2cm} \norm{w-\ell} & \geq \, \norm{ (\bar{B} x_1 + b_r x_{r,1}) - (\bar{B} x_2 + b_r x_{r,2}) }  = \norm{ \bar{B} (x_1 -  x_2) + b_r(x_{r,1} - x_{r,2}) }  \\ \nonumber
              & \geq \, \norm{b_r^*} \vert x_{r,1}-x_{r,2} \vert = \, \norm{b_r^*} \width(e_r,P)
\end{align}
holds, and so does \eref{erwidth}.

After branching on $e_r,\dots,e_{i+1}$, each subproblem is defined by a matrix formed of the first $i$ columns of $B$, and 
bound vectors $\ell_i$ and $w_i$, which are translates of $\ell$ and $w$ by the same vector. 
Hence the above proof implies that the width along $e_i$ in each of these subproblems is at most 
\beq \label{eiwidth} 
\norm{w-\ell}/\norm{b_i^*},
\eeq
and this completes the proof.
\qed

Our Lemma \ref{fractionlemma} builds on Furst and Kannan's Lemma 1 in \cite{FK89}, 
with inequality \eref{n} also being a direct generalization.
\ble \label{fractionlemma} 
For a positive integer $k$, let $\epsilon_R$ be the fraction of $A \in \gm$ with 
$\lambda_1(L_R(A)) \leq k, $ and $\epsilon_N$ be the fraction of $A \in \gmp$ with 
$\lambda_1(L_N(A)) \leq k. $ 
 Then 
\beq \label{r}
\epsilon_R \leq \dfrac{(2k+1)^{n+m}}{M^m},
\eeq
and 
\beq \label{n}
\epsilon_N \leq \dfrac{(2k+1)^n}{M^m \chi_{m,n}(M)}.
\eeq
\ele
\vspace{-0.4cm}
\pf{} We first prove \eref{n}. For $v, \,$ a fixed nonzero vector in $\zad{n}, \,$ consider the equation 
\beq \label{Av} 
Av = 0.
\eeq
There are at most $M^{m(n-1)}$ matrices in $\gmp$ that satisfy \eref{Av}:
if the components of $n-1$ columns of $A$ are fixed, then the components of the  column corresponding to a nonzero entry of $v$ 
are determined from \eref{Av}. 
The number of vectors $v$ in $\mathbb{Z}^{n}$ with $\norm{v} \leq k$ is at most $(2k+1)^n, $ and the 
number of matrices in $\gmp$ is  $M^{mn} \chi_{m,n}(M)$. Therefore 
$$
\eps_N \leq \dfrac{(2k+1)^nM^{m(n-1)}}{M^{mn} \chi_{m,n}(M)}=\dfrac{(2k+1)^n}{M^{m} \chi_{m,n}(M)}.
$$
For \eref{r}, note that $(v_1; v_2) \in \zad{m + n}$ is a nonzero vector in $L_R(A),$ iff $v_2 \neq 0, \,$ and 
\beq \label{Av2}
A v_2 = v_1.
\eeq
An argument like the one in the proof of \eref{n} shows that 
for fixed $(v_1; v_2) \in \zad{m + n}$ with $v_2 \neq 0, \,$ there are at most 
$M^{m(n-1)}$ matrices in $\gm$ that satisfy \eref{Av2}. 
The number of vectors in $\zad{n+m}$ with norm at most $k$ is at most $(2k+1)^{n+m}, $ so
$$\eps_R \leq \dfrac{(2k+1)^{n+m}M^{m(n-1)}}{M^{mn}}=\dfrac{(2k+1)^{n+m}}{M^{m}}.$$ 
\qed

\nin \pf{ of Theorems \ref{thm-rkz} and \ref{thm-lll}}
For part \eref{thm-rkz-r} in Theorem \ref{thm-rkz}, let $b_1^*, \dots, b_n^*$ be the Gram-Schmidt orthogonalization of the columns of $(A;I)U$. 
Lemma \ref{widthlemma} implies that the number of nodes generated by reverse B\&B applied to $Q_R$ is at most one, if 
\beq \label{bist}
\norm{b_i^*} > \norm{(w_1; w_2)  - (\ell_1; \ell_2)} 
\eeq
for $i=1, \dots, n$. Since the columns of $(A; I)U$ form an RKZ reduced basis of $L_R(A), \,$ \eref{bistci} implies
\beq \label{bilam} 
\norm{b_i^*} \geq \lambda_1(L_R(A))/{C_i}, 
\eeq
so \eref{bist} holds, when 
\beq \label{lambdacn}
\lambda_1(L_R(A)) > C_i \norm{(w_1; w_2)  - (\ell_1; \ell_2)} 
\eeq
does for $i=1, \dots, n$, which is in turn implied by 
\beq \label{lambdacn10}
\lambda_1(L_R(A)) > \gamma_n \norm{(w_1; w_2)  - (\ell_1; \ell_2)}. 
\eeq
By Lemma \ref{fractionlemma} \eref{lambdacn10} is true for all, but at most a fraction of $\eps_R$ of $A \in \gm$ if 
\beq
 M>\dfrac{( \lfloor 2\gamma_{n} \norm{ (w_1; w_2)-(\ell_1; \ell_2)}+1 \rfloor)^{(m+n)/m}}{\eps_R^{1/m}},
\eeq
and using the known estimate 
$\gamma_n \leq 1 + n/4 $ (see for instance \cite{M03}) , setting $\eps_R = 1/2^n, \,$ and doing some algebra yields the required result. 

The proof of part \eref{thm-rkz-n} of Theorem \ref{thm-rkz} is along the same lines: 
now \\ $b_1^*, \dots, b_{n-m}^*$ is the Gram-Schmidt orthogonalization of the columns of $B$, which is an RKZ reduced basis of $L_N(A)$. 
Lemma \ref{widthlemma}, and the reducedness of $B$ implies that the number of nodes generated by reverse B\&B applied to $Q_N$ is at most one, if 
\beq \label{lambdacn2}
\lambda_1(L_N(A)) > \gamma_{n-m} \norm{w_2 - \ell_2},
\eeq
\co{
Lemma \ref{widthlemma} implies that the number of nodes generated by reverse B\&B applied to $Q_N$ is at most one, if 
\beq \label{bist2}
\norm{b_i^*} > \norm{w_2 - \ell_2} 
\eeq
for $i=1, \dots, n-m$. By the reducedness of $B$,
\beq
\norm{b_i^*} \geq \lambda_1(L_N(A))/{C_i}, 
\eeq
so \eref{bist2} holds, when 
\beq \label{lambdacn2}
\lambda_1(L_N(A)) > C_{n-m} \norm{w_2 - \ell_2}
\eeq
does. }
and by Lemma \ref{fractionlemma} \eref{lambdacn2} is true for all, but at most a fraction of $\eps_N$ of $A \in \gm$ if 
\beq
 M>\dfrac{( \lfloor 2\gamma_{n-m} \norm{w_2-\ell_2}+1 \rfloor)^{n/m}}{\eps_N^{1/m} \chi_{m,n}(M)^{1/m}}.
\eeq
Then simple algebra, and using $\chi_{m,n}(M) \geq 1/2$ completes the proof. 

The proof of Theorem \ref{thm-lll} is an almost verbatim copy, now using the estimate 
\eref{bist2i} to lower bound $\norm{b_i^*}$. 
\qed

\co{and using the known estimate 
$C_n \leq 2n/3, $ setting $\eps_R = 1/2^n, \,$ and doing some algebra yields the required result. 
}


\pf{of Proposition \ref{widthcor} } 
To see \eref{widthrangerkz}, we start with 
\beq 
\width(e_n, P) \leq \norm{w-\ell}/\norm{b_n^*}
\eeq
from \eref{erwidth}, combine it with the lower bound on $\norm{b_n^*}$ from \eref{brrkz}, and the fact that 
\begin{align}  \label{detlra}
\mydet L_R(A) & = \mydet (A A^T+I)^{1/2}, 
\end{align}
which follows from the definition of $L_R(A)$. The proof of \eref{widthnullrkz} is analogous, but now we need to use 
\begin{align}  \label{detlna}
\mydet L_R(A) & = \mydet (AA^T)^{1/2}/\gcd(A),
\end{align}
whose proof can be found in \cite{Cassels} for instance. To prove the claims about the LLL-reformulations, we need to use 
\eref{brlll} in place of \eref{brrkz}. 
\qed

\pf{of Proposition \ref{actual}}
\vspace{-0.3cm}
Let $N(n,k)$ denote the number of integral points in the $n$-dimensional ball of radius $k$. 
In the previous proofs we used $(2k+1)^n$ as an upper bound for $N(n,k)$. The proof of Part \eref{thm-rkz-n} of Theorem \ref{thm-rkz} 
actually implies that 
when 
\beq \label{MCnm} 
 M>\dfrac{( N(n, \lceil \gamma_{n-m} \norm{w_2-\ell_2} \rceil)^{1/m}}{\eps_N^{1/m} \chi_{m,n}(M)},
\eeq
then for all, but at most a fraction 
of $\eps_N$ of $A \in \gm\,$  reverse B\&B solves the nullspace reformulation of \eref{Ab} at the rootnode. 

We use $\chi_{m,n}(M) \geq 1/2, \,$ Blichfeldt's upper bound 
\begin{equation} \label{Bl}
C_i \leq \dfrac{2}{\pi} {\Gamma \left( {\dfrac{i+4}{2}} \right) }^{2/i},
\end{equation}
from \cite{Bl14} to bound $\gamma_{n-m}$ in \eref{MCnm}, 
dynamic programming to exactly find the values of $N(n,k),$ and the values 
$\eps_N = 0.1, \,$ and $\eps_N = 0.01 \,$ to obtain Table \ref{mnMtable90}. 

We note that in general $N(n,k)$ is hard to compute, or find good upper bounds for; however for small values of $n$ and $k$ 
a simple dynamic programming algorithm finds the exact value quickly.\qed

\co{
First we note that if $b_1^*, \dots, b_r^*$ is an RKZ reduced, or LLL reduced basis of the lattice $L$, then by \cite{LLS90} 
\beq
\norm{b_r^*} \geq \dfrac{(\mydet L)^{1/r}}{\sqrt{r}};
\eeq
if it is an LLL reduced basis, then multiplying the inequalities 
\beq
\norm{b_i^*} \leq 2^{(r-i)/2} \norm{b_r^*} \; (i=1, \dots, r),
\eeq
and using $\norm{b_1^*} \dots \norm{b_r^*} = \mydet L$ gives 
\beq
\norm{b_r^*} \geq \dfrac{(\mydet L)^{1/r}}{2^{(r-1)/4}}. 
\eeq
}

\nin{\bf Acknowledgement} We thank Van Vu for pointing out reference \cite{BVW09}.

\bibliographystyle{plain}
\bibliography{C:/bibfiles/IP_Refs}

\begin{thebibliography}{10}

\bibitem{ABHLS00}
Karen Aardal, Robert~E. Bixby, Cor A.~J. Hurkens, Arjen~K. Lenstra, and Job~W.
  Smeltink.
\newblock Market split and basis reduction: Towards a solution of the
  {C}ornu{\'{e}}jols-{D}awande instances.
\newblock {\em INFORMS Journal on Computing}, 12(3):192--202, 2000.

\bibitem{AHL00}
Karen Aardal, Cor A.~J. Hurkens, and Arjen~K. Lenstra.
\newblock Solving a system of linear {D}iophantine equations with lower and
  upper bounds on the variables.
\newblock {\em Mathematics of Operations Research}, 25(3):427--442, 2000.

\bibitem{AL04}
Karen Aardal and Arjen~K. Lenstra.
\newblock Hard equality constrained integer knapsacks.
\newblock {\em Mathematics of Operations Research}, 29(3):724--738, 2004.

\bibitem{AL06}
Karen Aardal and Arjen~K. Lenstra.
\newblock Erratum to: Hard equality constrained integer knapsacks.
\newblock {\em Mathematics of Operations Research}, 31(4):846, 2006.

\bibitem{Bl14}
Hans~Frederik Blichfeldt.
\newblock A new principle in the geometry of numbers, with some applications.
\newblock {\em Transactions of the American Mathematical Society},
  15(3):227--235, 1914.

\bibitem{BVW09}
Jean Bourgain, Van~H. Vu, and Philip~Matchett Wood.
\newblock On the singularity probability of discrete random matrices.
\newblock {\em Journal of Functional Analysis}, to appear.

\bibitem{Cassels}
J.~W.~S. Cassels.
\newblock {\em An introduction to the geometry of numbers}.
\newblock Springer, 1997.

\bibitem{C80}
Va{\v{s}}ek Chv{\'{a}}tal.
\newblock Hard knapsack problems.
\newblock {\em Operations Research}, 28(6):1402--1411, 1980.

\bibitem{CD99}
G{\'{e}}rard Cornu{\'{e}}jols and Milind Dawande.
\newblock A class of hard small 0--1 programs.
\newblock {\em INFORMS Journal on Computing}, 11(2):205--210, 1999.

\bibitem{CKL06}
G{\'{e}}rard Cornu{\'{e}}jols, Miroslav Karamanov, and Yanjun Li.
\newblock Early estimates of the size of branch-and-bound trees.
\newblock {\em INFORMS Journal on Computing}, 18(1):86--96, 2006.

\bibitem{F86}
Alan Frieze.
\newblock On the \mbox{Lagarias-Odlyzko} algorithm for the subset sum problem.
\newblock {\em SIAM Journal on Computing}, 15:536--540, 1986.

\bibitem{FK89}
Merrick Furst and Ravi Kannan.
\newblock Succinct certificates for almost all subset sum problems.
\newblock {\em SIAM Journal on Computing}, 18:550 -- 558, 1989.

\bibitem{GNS98Complex}
Zonghao Gu, George~L. Nemhauser, and Martin W.~P. Savelsbergh.
\newblock Lifted cover inequalities for 0--1 integer programs: Complexity.
\newblock {\em INFORMS J. on Computing}, 11:117--123, 1998.

\bibitem{J74}
Robert~G. Jeroslow.
\newblock Trivial integer programs unsolvable by \bnb.
\newblock {\em Mathematical Programming}, 6:105--109, 1974.

\bibitem{K87}
Ravi Kannan.
\newblock Minkowski's convex body theorem and integer programming.
\newblock {\em Mathematics of Operations Research}, 12(3):415--440, 1987.

\bibitem{KZ1873}
A.~Korkine and G.~Zolotarev.
\newblock Sur les formes quadratiques.
\newblock {\em Mathematische Annalen}, 6:366--389, 1873.

\bibitem{KP09}
Bala Krishnamoorthy and G{\'{a}}bor Pataki.
\newblock Column basis reduction and decomposable knapsack problems.
\newblock {\em Discrete Optimization}, 6:242--270, 2009.

\bibitem{LLS90}
Jeffrey~C. Lagarias, Hendrik~W. Lenstra, and Claus~P. Schnorr.
\newblock Korkine-{Z}olotarev bases and successive minina of a lattice and its
  reciprocal lattice.
\newblock {\em Combinatorica}, 10(4):333--348, 1990.

\bibitem{LO85}
Jeffrey~C. Lagarias and Andrew~M. Odlyzko.
\newblock Solving low-density subset sum problems.
\newblock {\em Journal of ACM}, 32:229--246, 1985.

\bibitem{LD60}
A.~H. Land and Alison~G. Doig.
\newblock An automatic method for solving discrete programming problems.
\newblock {\em Econometrica}, 28:497--520, 1960.

\bibitem{LLL82}
Arjen~K. Lenstra, Hendrik~W. {Lenstra, Jr.}, and L{\'{a}}szl{\'{o}}
  Lov{\'{a}}sz.
\newblock Factoring polynomials with rational coefficients.
\newblock {\em Mathematische Annalen}, 261:515--534, 1982.

\bibitem{L83}
Hendrik~W. {Lenstra, Jr.}
\newblock Integer programming with a fixed number of variables.
\newblock {\em Mathematics of Operations Research}, 8:538--548, 1983.

\bibitem{LS92}
L{\'{a}}szl{\'{o}} Lov{\'{a}}sz and Herbert~E. Scarf.
\newblock The generalized basis reduction algorithm.
\newblock {\em Mathematics of Operations Research}, 17:751--764, 1992.

\bibitem{M03}
Jacques Martinet.
\newblock {\em Perfect Lattices in Euclidean Spaces}.
\newblock Springer-Verlag, Berlin, 2003.

\bibitem{ML04}
Sanjay Mehrotra and Zhifeng Li.
\newblock On generalized branching methods for mixed integer programming.
\newblock {\em Research Report, Department of Industrial Engineering,
  Northwestern University}, 2004.

\bibitem{NTL}
Victor Shoup.
\newblock {NTL}: {A} {N}umber {T}heory {L}ibrary, 1990.
\newblock http://www.shoup.net.

\end{thebibliography}
\end{document}